\documentclass[12pt]{amsart}

\usepackage{color,fancyhdr}
\usepackage[all]{xy}

\usepackage{amsmath,amsthm,mathrsfs,amsfonts,
amssymb,times,latexsym,mathabx}

\usepackage[usenames,dvipsnames,svgnames,table]{xcolor}

\DeclareMathAlphabet{\curly}{U}{rsfs}{m}{n}  

\newtheorem{theorem}{Theorem}[section]

\theoremstyle{definition}

\newtheorem{corollary}[theorem]{Corollary}

\theoremstyle{problem}
\newtheorem{problem}[theorem]{Problem}
\usepackage{xcolor}

\numberwithin{equation}{section}

\makeatletter
\renewcommand{\pmod}[1]{\allowbreak\mkern7mu({\operator@font mod}\,\,#1)}
\makeatother

\newcommand{\be}{\begin{equation}}
\newcommand{\ee}{\end{equation}}

\renewcommand{\le}{\leqslant}

\renewcommand{\ge}{\geqslant}

\begin{document}

\title[On function $SX$ of additive complements]
{On function $SX$ of additive complements}

\author{Jin-Hui Fang}
\address{Department of Mathematics, Nanjing University of Information Science $\&$ Technology, Nanjing 210044, PR China}
\email{fangjinhui1114@163.com}

\author{Csaba S\'{a}ndor*}
\address{Institute of Mathematics, Budapest University of Technology and Economics, Egry J\'ozsef utca 1, 1111 Budapest, Hungary; Department of Computer Science and Information Theory, Budapest University of Technology and Economics, M\H{u}egyetem rkp. 3., H-1111 Budapest, Hungary; MTA-BME Lend\"{u}let Arithmetic Combinatorics Research Group, ELKH, M\H{u}egyetem rkp. 3., H-1111 Budapest, Hungary}
\email{csandor@math.bme.hu}
\thanks{* Corresponding author.}
\thanks{The first author is supported by the National Natural Science Foundation of China, Grant No. 12171246 and the Natural Science Foundation of Jiangsu Province, Grant No. BK20211282. The second author is supported by the NKFIH Grants No. K129335.}
\keywords{Additive complements, Perfect, Infimum}
\subjclass[2010]{Primary 11B13, Secondary 11B34}
\date{\today}%

\begin{abstract}
Two sets $A,B$ of nonnegative integers are called \emph{additive complements}, if all sufficiently large integers can be expressed as the sum of two elements from $A$ and $B$. We further call $A,B$ \emph{perfect additive complements} if every nonnegative integer can be uniquely expressed as the sum of two elements from $A$ and $B$. Let $A(x)$ be the counting function of $A$. In this paper, we focus on the function $SX$, where $SX=\limsup_{x\rightarrow\infty}\frac{\max\{A(x),B(x)\}}{\sqrt{x}}$ was introduced by Erd\H os and Freud in 1984. As a main result, we determine the value of $SX$ for perfect additive complements and further fix the infimum. We also give the absolute lower bound of $SX$ for additive complements.
\end{abstract}
\maketitle

\section{\bf Introduction}

Two sets $A,B$ of nonnegative integers are called \emph{additive complements}, if all sufficiently large integers can be expressed as the sum of two elements from $A$ and $B$. Let $A(x)$ be the counting function of $A$. In 1959, Narkiewicz proved the following remarkable result.

\noindent\textbf{Theorem A.}
For additive complements $A$ and $B$ with $A(x)B(x)=(1+o(1))x$ we have $$\lim_{x\to \infty}\frac{\log \min \{ A(x),B(x)\}}{\log x}=0.$$ \vskip2mm

 We firstly generalize Narkiewicz' result on additive complements.
\vskip2mm

\begin{theorem}\label{thm1} For additive complements $A$ and $B$ with $\limsup_{x\to \infty}\frac{A(x)B(x)}{x}= 1+\delta$, where
\begin{eqnarray*}0\le \delta <C_0
=\frac{1}{2}(-3-\sqrt{2}+\sqrt{3+12\sqrt{2}})=0.027315\cdots,\end{eqnarray*}
we have \begin{eqnarray*}\limsup_{x\to \infty}\frac{\log \min \{ A(x),B(x)\}}{\log x}\le f(\delta),\end{eqnarray*}
where $f(\delta)$ is defined as
\begin{eqnarray*}
f(\delta)=\log_2 \frac{(2\delta +2)^2}{3-2\delta +\sqrt{(3-2\delta)^2-(2\delta +2)^3}}=f(\delta).\end{eqnarray*}
\end{theorem}
\vskip2mm

\noindent\textbf{{Remark.}} We note that $f(0)=0$, $f(C_0)=0.5$ and $f(\delta)$ is strictly increasing for $0\le \delta \le C_0$.
\vskip2mm

As a direct consequence we get the following.

\begin{corollary}
For additive complements $A$ and $B$ with $\limsup_{x\to \infty}\frac{A(x)B(x)}{x}<1+C_0$ we have $\lim_{x\to \infty }\frac{\min \{ A(x),B(x)\} }{\sqrt{x}}=0$.
\end{corollary}

For any additive complements $A$ and $B$ we have $A(x)B(x)\ge x-c_1$. Hence

\begin{corollary}\label{140}
For additive complements $A$ and $B$ with $\limsup_{x\to \infty}\frac{A(x)B(x)}{x}<1+C_0$ we have $\lim_{x\to \infty }\frac{\max \{ A(x),B(x)\} }{\sqrt{x}}=\infty $.
\end{corollary}

For additive complements $A$ and $B$, let \begin{eqnarray*}SX(A,B)=\limsup_{x\rightarrow\infty}\frac{\max\{A(x),B(x)\}}{\sqrt{x}},\end{eqnarray*}
which was introduced by Erd\H os and Freud in \cite{Erdos}. In this paper, we mainly focus our attention on $SX(A,B)$. \vskip2mm

\begin{theorem}\label{thm2} For additive complements $A$ and $B$, we have
\begin{eqnarray*}SX(A,B)\ge \sqrt{1+C_0}=1.013565....\end{eqnarray*}
\end{theorem}\vskip2mm

Denote the sets $A,B$ by \emph{perfect additive complements} if every nonnegative integer can be uniquely expressed as the sum of two elements from $A$ and $B$. In \cite{Sandor} we fixed the structure of \emph{perfect additive complements} as follows:
\begin{eqnarray}\label{10031}
&&A=\{\epsilon_0+\epsilon_2m_1m_2+\cdots+\epsilon_{2k-2}m_1\cdots m_{2k-2}+\cdots, \epsilon_{2i}=0,1,\cdots,m_{2i+1}-1\},\nonumber\\
&&B=\{\epsilon_1m_1+\epsilon_3m_1m_2m_3+\cdots+\epsilon_{2k-1}m_1\cdots m_{2k-1}+\cdots, \epsilon_{2i-1}=0,1,\cdots,m_{2i}-1\},
\end{eqnarray}(or $A,B$ interchanged), where $m_1,m_2,\cdots$ are integers no less than two.
We determine the value of $SX(A,B)$ for \emph{perfect additive complements} with the form \eqref{10031}, that is:\vskip2mm
\begin{theorem}\label{thm3}
\begin{eqnarray*}
SX(A,B)&=& \limsup_{s\to \infty}\max \{ \frac{m_1m_3\cdots m_{2s-1}}
{\sqrt{(m_1-1)+(m_3-1)m_1m_2+\cdots+(m_{2s-1}-1)m_1m_2\cdots m_{2s-2}}},\\
&&\frac{m_2m_4\cdots m_{2s}}
{\sqrt{(m_2-1)m_1+(m_4-1)m_1m_2m_3+\cdots+(m_{2s}-1)m_1m_2\cdots m_{2s-1}}}\}.
\end{eqnarray*}
\end{theorem}\vskip2mm

We further fix the infimum of $SX$ for perfect additive complements.\vskip2mm

\begin{theorem}\label{thm4}
The infimum of $SX(A,B)$ for perfect additive complements is $\sqrt[4]{4.5}$.
\end{theorem}
\vskip2mm

\noindent\textbf{Remark.} We find that the supremum of $SX$ for perfect additive complements is $\infty$ by taking $m_{2i-1}=2$ and $m_{2i}=3$ for all positive integer $i$.
\vskip3mm

We pose two problems for further research.

\begin{problem}
Is it true that for any $\delta>0$ there are additive complements $A$ and $B$ such that $\limsup_{x\to \infty}\frac{A(x)B(x)}{x}\le 1+\delta $ and $\liminf_{x\to \infty }\frac{\log \min \{ A(x),B(x)\} }{\log x}>0$?
\end{problem}

\begin{problem}
Is it true that for any additive complements $A$ and $B$ we have $SX(A,B)\ge \sqrt[4]{4.5}$?
\end{problem}
\section{Proof of Main Results}\label{sec:proofofMT}

\noindent\textbf{Proof of Theorem \ref{thm1}.} 
Suppose that $\limsup_{x\to \infty}\frac{A(x)B(x)}{x}=1+\delta$, where $0\le \delta < C_0$. Then
\begin{eqnarray*}
(1+\delta +o(1))x&\ge& \sum_{(a,b), a\in A, b\in B,a\le x,b\le x}1=\sum_{(a,b), a\in A, b\in B,a\le x,b\le x, a+b\le x}1+\sum_{(a,b), a\in A, b\in B,a\le x,b\le x, a+b>x}1\\
&\ge& x+(A(x)-A(\frac{x}{2}))(B(x)-B(\frac{x}{2})).
\end{eqnarray*}
Hence
\begin{eqnarray*}
(\delta+o(1))x&\ge&(A(x)-A(\frac{x}{2}))(B(x)-B(\frac{x}{2}))
=A(x)B(x)+A(\frac{x}{2})B(\frac{x}{2})-A(x)B(\frac{x}{2})-B(x)A(\frac{x}{2})\\
&\ge& x+0.5x-A(x)B(\frac{x}{2})-B(x)A(\frac{x}{2}),
\end{eqnarray*}
that is
$$
\frac{A(\frac{x}{2})}{A(x)}+\frac{B(\frac{x}{2})}{B(x)}=\frac{A(x)B(\frac{x}{2})+B(x)A(\frac{x}{2})}{A(x)B(x)}\ge \frac{(1.5-\delta +o(1))x}{(1+\delta +o(1))x}= \frac{1.5-\delta}{1+\delta}+o(1).
$$
It follows from
\begin{eqnarray*}
\frac{B(\frac{x}{2})}{B(x)}\le \frac{\frac{(1+\delta +o(1))\frac{x}{2}}{A(\frac{x}{2})}}{\frac{x}{A(x)}}=\frac{A(x)}{A(\frac{x}{2})}\left(\frac{1+\delta }{2}+o(1)\right)
\end{eqnarray*}
that
\begin{eqnarray*}
\frac{1.5-\delta}{1+\delta}+o(1)\le \frac{A(\frac{x}{2})}{A(x)}+\frac{A(x)}{A(\frac{x}{2})}\left(\frac{1+\delta }{2}+o(1)\right).
\end{eqnarray*}
Thus,
\begin{eqnarray*}0\le 2(1+\delta)+o(1)-(3-2\delta +o(1))\frac{A(x)}{A(\frac{x}{2})}+((1+\delta )^2+o(1))\left(\frac{A(x)}{A(\frac{x}{2})}\right)^2.
\end{eqnarray*}
It is easy to check that the quadratic polynomial $p_{\delta}(z)=2(1+\delta)+(3-2\delta)z+(1+\delta)^2z^2$ is a perfect square for $\delta _0 =0.027357\dots $ and there are two real roots $$r_1(\delta)=\frac{3-2\delta-\sqrt{(3-2\delta)^2-8(1+\delta)^3}}{2(1+\delta)^2}<r_2(\delta)=\frac{3-2\delta+\sqrt{(3-2\delta)^2-8(1+\delta)^3}}{2(1+\delta)^2}$$ for $0\le \delta <\delta _0$. A simple calculation gives that $r_2(C_0)=\sqrt{2}=1.41421\dots $ and $r_1(C_0)=1.37661\dots $. The function $r_1(x)$ is monotone increasing in $[0,\delta _0[$ and the function $r_2(x)$ is monotone decreasing in $[0,\delta _0[$. Since \begin{eqnarray*}\frac{A(n+1)}{A(\frac{n+1}{2})}-\frac{A(n)}{A(\frac{n}{2})}\to 0,\end{eqnarray*} there exists a positive number $x_0$ such that
\begin{eqnarray*}
\frac{A(x)}{A(\frac{x}{2})}\ge (r_2(\delta )-o(1))^x \hskip3mm  \mbox{or} \hskip3mm \frac{A(x)}{A(\frac{x}{2})}\le (r_1(\delta)-o(1))^x \hskip2mm  \mbox{for} \hskip2mm  x\ge x_0.\end{eqnarray*} In the first case $A(2^n)\gg (r_2(\delta)-o(1))^n$, so $\liminf _{x\to \infty }\frac{\log A(x)}{\log x}\ge \log _2 r_2(\delta)$. Then $\limsup_{x\to \infty }\frac{\log B(x)}{\log x}\le 1-\log _2 r_2(\delta)$. In the second case $A(2^n)\ll (r_1(\delta)-o(1))^n$, so $\limsup_{x\to \infty }\frac{\log A(x)}{\log x}\le \log _2 r_1(\delta)$. To sum up,
$$\limsup _{x\to \infty}\frac{\min \{ A(x),B(x)\}}{\log x}\le \max \{ \log _2 r_1(\delta),1-\log _2 r_2(\delta)\}.$$
Clearly,
$$r_2(\delta)r_1(\delta)=\frac{2}{1+\delta}.$$
Thus,
$$\max \{ \log _2 r_1(\delta),1-\log _2 r_2(\delta)\}=1-\log _2 r_2(\delta)=f(\delta).$$
\vskip2mm

This completes the proof of Theorem \ref{thm1}. \hfill$\Box$\\

\noindent\textbf{Proof of Theorem \ref{thm2}.} If $\limsup_{x\to \infty}\frac{A(x)B(x)}{x}<1+C_0$, then by Corollary \ref{140} we have $SX(A,B)=\infty$. If $\limsup_{x\to \infty}\frac{A(x)B(x)}{x}\ge 1+C_0$, then
\begin{eqnarray*}
SX(A,B)=\limsup_{x\rightarrow\infty}\frac{\max\{A(x),B(x)\}}{\sqrt{x}}
\ge \sqrt{\limsup_{x\to \infty}\frac{A(x)B(x)}{x}}\ge \sqrt{1+C_0}=1.013\dots .\end{eqnarray*}

This completes the proof of Theorem \ref{thm2}. \hfill$\Box$\\

\noindent\textbf{Proof of Theorem \ref{thm3}.}
If $$x=(m_1-1)+(m_3-1)m_1m_2+\cdots+(m_{2s-1}-1)m_1m_2\cdots m_{2s-2},$$
then
$$A(x)=m_1m_3\cdots m_{2s-1}$$
and if
$$x=(m_2-1)m_1+(m_4-1)m_1m_2m_3+\cdots+(m_{2s}-1)m_1m_2\cdots m_{2s-1}$$
then
$$B(x)=m_2m_4\cdots m_{2s}.$$
Hence
\begin{eqnarray*}
SX(A,B)&\ge& \limsup _{s\to \infty}\max \{ \frac{m_1m_3\cdots m_{2s-1}}
{\sqrt{(m_1-1)+(m_3-1)m_1m_2+\cdots+(m_{2s-1}-1)m_1m_2\cdots m_{2s-2}}},\\
&&\frac{m_2m_4\cdots m_{2s}}
{\sqrt{(m_2-1)m_1+(m_4-1)m_1m_2m_3+\cdots+(m_{2s}-1)m_1m_2\cdots m_{2s-1}}}\}.
\end{eqnarray*}
Clearly,
\begin{eqnarray*}SX(A,B)=\limsup_{x\rightarrow\infty}\frac{\max\{A(x),B(x)\}}{\sqrt{x}}=\max \{ \limsup_{x\to \infty }\frac{A(x)}{\sqrt{x}},\limsup_{x\to \infty }\frac{B(x)}{\sqrt{x}} \}. \end{eqnarray*}
To finish the proof it is enough to show that for
\begin{eqnarray}\label{10032}y&=&(m_1-1)+(m_3-1)m_1m_2+\dots+(m_{2j-1}-1)m_1m_2\cdots m_{2j-2}+\nonumber\\
&&\delta_{2j+1}m_1m_2\cdots m_{2j}+\delta_{2j+3}m_1m_2\cdots m_{2j+2}+\cdots+\delta_{2k+1}m_1m_2\cdots m_{2k},
\end{eqnarray}
where $0\le \delta _{2j+1}<m_{2j+1}-1$ and $1\le \delta _{2k+1}\le m_{2k+1}-1$ we have
$$\frac{A(y)}{\sqrt{y}}\le \frac{A(y+m_1m_2\cdots m_{2j})}{\sqrt{y+m_1m_2\cdots m_{2j}}}$$
and for
\begin{eqnarray}\label{10033}z&=&(m_2-1)m_1+(m_4-1)m_1m_2m_3+\dots +(m_{2j}-1)m_1m_2\cdots m_{2j-1}+\nonumber\\
&&\delta _{2j+2}m_1m_2\cdots m_{2j+1}+\delta_{2j+4}m_1m_2\cdots m_{2j+3}+\cdots+\delta_{2k}m_1m_2\cdots m_{2k-1},
\end{eqnarray}
where $0\le \delta _{2j+2}<m_{2j+2}-1$ and $1\le \delta _{2k}\le m_{2k}-1$ we have
$$\frac{A(z)}{\sqrt{z}}\le \frac{A(z+m_1m_2\cdots m_{2j+1})}{\sqrt{z+m_1m_2\cdots m_{2j+1}}}.$$
For $y$ with the form \eqref{10032}, we have $$A(y+m_1m_2\cdots m_{2j})=A(y)+m_1m_3\dots m_{2j-1}.$$ So we need to show that
$$\frac{A(y)+m_1m_3\dots m_{2j-1}}{A(y)}\ge \sqrt{\frac{y+m_1m_2\cdots m_{2j}}{y}}.$$
Since
$$\sqrt{\frac{y+m_1m_2\cdots m_{2j}}{y}}\le 1+\frac{m_1m_2\cdots m_{2j}}{2y},$$
it suffices to prove that
$$\frac{m_1m_3\dots m_{2j-1}}{A(y)}\ge \frac{m_1m_2\cdots m_{2j}}{2y},$$
that is
$$2y\ge A(y)m_2m_4\dots m_{2j},$$
which follows from $y\ge \delta _{2k+1}m_1m_2\cdots m_{2k}$ and $A(y)\le (\delta _{2k+1}+1)m_1m_3\dots m_{2k-1}$.

For $z$ with the form \eqref{10033}, the proof is similar.\vskip2mm

This completes the proof of Theorem \ref{thm3}. \hfill$\Box$\\

\noindent\textbf{Proof of Theorem \ref{thm4}.} Let $m_1^{(k)}, m_2^{(k)}$ be positive integers with
$\frac{m_2^{(k)}}{m_1^{(k)}}\rightarrow \sqrt{2}$ as $k\to \infty$ and $m_n^{(k)}=2$ for all $k\ge 1$ and $n\ge 3$. It is easy to check that
$$\lim_{k\to \infty} SX(A^{(k)},B^{(k)})=\sqrt[4]{4.5}.$$ By Theorem \ref{thm3}, it suffices to prove that
\begin{eqnarray*}
\min \{ \liminf_{s\to \infty }&& ( \frac{1}{m_{2s}}(1-\frac{1}{m_{2s-1}}) \frac{m_{2s}}{m_{2s-1}}\frac{m_{2s-2}}{m_{2s-3}}\cdots \frac{m_2}{m_1}+\\
&&\frac{1}{m_{2s-1}^2}\frac{1}{m_{2s-2}}(1-\frac{1}{m_{2s-3}}) \frac{m_{2s-2}}{m_{2s-3}}\frac{m_{2s-4}}{m_{2s-5}}\cdots \frac{m_2}{m_1}+\\
&&\frac{1}{m_{2s-1}^2}\frac{1}{m_{2s-3}^2}\frac{1}{m_{2s-4}}(1-\frac{1}{m_{2s-5}}) \frac{m_{2s-4}}{m_{2s-5}}\frac{m_{2s-6}}{m_{2s-7}}\cdots \frac{m_2}{m_1}+\\
&&\frac{1}{m_{2s-1}^2}\frac{1}{m_{2s-3}^2}\frac{1}{m_{2s-5}^2}\frac{1}{m_{2s-6}}
(1-\frac{1}{m_{2s-7}})\frac{m_{2s-6}}{m_{2s-7}}\frac{m_{2s-8}}{m_{2s-9}}\cdots \frac{m_2}{m_1}+\cdots +\\
&&\frac{1}{m_{2s-1}^2}\frac{1}{m_{2s-3}^2}\cdots \frac{1}{m_{2s-(2k-1)}^2}\frac{1}{m_{2s-2k}}
(1-\frac{1}{m_{2s-(2k+1)}})\frac{m_{2s-2k}}{m_{2s-(2k+1)}}\frac{m_{2s-(2k+2)}}{m_{2s-(2k+3)}}\cdots \frac{m_2}{m_1}+\cdots ),\\
\liminf_{s\to \infty}&&((1-\frac{1}{m_{2s}}) \frac{m_{2s-1}}{m_{2s}}\frac{m_{2s-3}}{m_{2s-2}}\cdots \frac{m_1}{m_2}+\\
&&\frac{1}{m_{2s}^2}(1-\frac{1}{m_{2s-2}}) \frac{m_{2s-3}}{m_{2s-2}}\frac{m_{2s-5}}{m_{2s-4}}\cdots \frac{m_1}{m_2}\\
&&+\frac{1}{m_{2s}^2}\frac{1}{m_{2s-2}^2}(1-\frac{1}{m_{2s-4}}) \frac{m_{2s-5}}{m_{2s-4}}\frac{m_{2s-7}}{m_{2s-6}}\cdots \frac{m_1}{m_2}\\
&&+\frac{1}{m_{2s}^2}\frac{1}{m_{2s-2}^2}\frac{1}{m_{2s-4}^2}(1-\frac{1}{m_{2s-6}}) \frac{m_{2s-7}}{m_{2s-6}}\frac{m_{2s-9}}{m_{2s-8}}\cdots \frac{m_1}{m_2}+\cdots +\\
&& \frac{1}{m_{2s}^2}\frac{1}{m_{2s-2}^2}\cdots \frac{1}{m_{2s-(2k-2)}^2}(1-\frac{1}{m_{2s-2k}}) \frac{m_{2s-(2k+1)}}{m_{2s-2k}}\frac{m_{2s-(2k+3)}}{m_{2s-(2k+2)}}\cdots \frac{m_1}{m_2} +\cdots ) \}\\
&&\le \frac{\sqrt{2}}{3}.
\end{eqnarray*}

Suppose that $m_{2s}\ge 3$ for infinitely many $s$. Let $m_{2s}\ge 3$. Then
\begin{eqnarray*}
&& ( \frac{1}{m_{2s}}(1-\frac{1}{m_{2s-1}}) \frac{m_{2s}}{m_{2s-1}}\frac{m_{2s-2}}{m_{2s-3}}\cdots \frac{m_2}{m_1}+\frac{1}{m_{2s-1}^2}\frac{1}{m_{2s-2}}(1-\frac{1}{m_{2s-3}}) \frac{m_{2s-2}}{m_{2s-3}}\frac{m_{2s-4}}{m_{2s-5}}\cdots \frac{m_2}{m_1}\\
&&+\frac{1}{m_{2s-1}^2}\frac{1}{m_{2s-3}^2}\frac{1}{m_{2s-4}}(1-\frac{1}{m_{2s-5}}) \frac{m_{2s-4}}{m_{2s-5}}\frac{m_{2s-6}}{m_{2s-7}}\cdots \frac{m_2}{m_1}\\
&&+\frac{1}{m_{2s-1}^2}\frac{1}{m_{2s-3}^2}\frac{1}{m_{2s-5}^2}\frac{1}{m_{2s-6}}
(1-\frac{1}{m_{2s-7}})\frac{m_{2s-6}}{m_{2s-7}}\frac{m_{2s-8}}{m_{2s-9}}\cdots \frac{m_2}{m_1}+\cdots +\\
&& \frac{1}{m_{2s-1}^2}\frac{1}{m_{2s-3}^2}\cdots \frac{1}{m_{2s-(2k-1)}^2}\frac{1}{m_{2s-2k}}
(1-\frac{1}{m_{2s-(2k+1)}})\frac{m_{2s-2k}}{m_{2s-(2k+1)}}\frac{m_{2s-(2k+2)}}{m_{2s-(2k+3)}}\cdots \frac{m_2}{m_1}+\cdots)\times \\
&&((1-\frac{1}{m_{2s}}) \frac{m_{2s-1}}{m_{2s}}\frac{m_{2s-3}}{m_{2s-2}}\cdots
\frac{m_1}{m_2}+\frac{1}{m_{2s}^2}(1-\frac{1}{m_{2s-2}}) \frac{m_{2s-3}}{m_{2s-2}}\frac{m_{2s-5}}{m_{2s-4}}\cdots \frac{m_1}{m_2}\\
&&+\frac{1}{m_{2s}^2}\frac{1}{m_{2s-2}^2}(  1-\frac{1}{m_{2s-4}}) \frac{m_{2s-5}}{m_{2s-4}}\frac{m_{2s-7}}{m_{2s-6}}\cdots \frac{m_1}{m_2}\\
&&+\frac{1}{m_{2s}^2}\frac{1}{m_{2s-2}^2}\frac{1}{m_{2s-4}^2}(1-\frac{1}{m_{2s-6}}) \frac{m_{2s-7}}{m_{2s-6}}\frac{m_{2s-9}}{m_{2s-8}}\cdots \frac{m_1}{m_2}+\cdots +\\
&&\frac{1}{m_{2s}^2}\frac{1}{m_{2s-2}^2}\cdots \frac{1}{m_{2s-(2k-2)}^2}(1-\frac{1}{m_{2s-2k}}) \frac{m_{2s-(2k+1)}}{m_{2s-2k}}\frac{m_{2s-(2k+3)}}{m_{2s-(2k+2)}}\cdots \frac{m_1}{m_2} +\cdots )=\\
&&(\frac{1}{m_{2s}}(1-\frac{1}{m_{2s}})(1-\frac{1}{m_{2s-1}}))+\\
&&(\frac{1}{m_{2s}}\frac{1}{m_{2s-1}}\frac{1}{m_{2s-2}}(1-\frac{1}{m_{2s}})
(1-\frac{1}{m_{2s-3}})+\frac{1}{m_{2s}^2}\frac{1}{m_{2s-1}}
(1-\frac{1}{m_{2s-1}})(1-\frac{1}{m_{2s-2}}))+ \\
&&(\frac{1}{m_{2s}}\frac{1}{m_{2s-1}}\frac{1}{m_{2s-2}}\frac{1}{m_{2s-3}}
\frac{1}{m_{2s-4}}(1-\frac{1}{m_{2s}})(1-\frac{1}{m_{2s-5}})\\
&&+\frac{1}{m_{2s}^2}\frac{1}{m_{2s-1}}\frac{1}{m_{2s-2}}\frac{1}{m_{2s-3}}(1-\frac{1}{m_{2s-1}})
(1-\frac{1}{m_{2s-4}})+\cdots \\
&&+\frac{1}{m_{2s}^2}\frac{1}{m_{2s-1}^2}\frac{1}{m_{2s-2}}(1-\frac{1}{m_{2s-2}})
(1-\frac{1}{m_{2s-3}}))\\
&&+(\frac{1}{m_{2s}}\frac{1}{m_{2s-1}}\frac{1}{m_{2s-2}}\frac{1}{m_{2s-3}}
\frac{1}{m_{2s-4}}\frac{1}{m_{2s-5}}\frac{1}{m_{2s-6}}(1-\frac{1}{m_{2s}})(1-\frac{1}{m_{2s-7}})\\
&&+\frac{1}{m_{2s}^2}\frac{1}{m_{2s-1}}\frac{1}{m_{2s-2}}\frac{1}{m_{2s-3}}
\frac{1}{m_{2s-4}}\frac{1}{m_{2s-5}}(1-\frac{1}{m_{2s-1}})(1-\frac{1}{m_{2s-6}})\\
&&+\frac{1}{m_{2s}^2}\frac{1}{m_{2s-1}^2}\frac{1}{m_{2s-2}}
\frac{1}{m_{2s-3}}\frac{1}{m_{2s-4}}(1-\frac{1}{m_{2s-2}})(1-\frac{1}{m_{2s-5}})\\
&&+\frac{1}{m_{2s}^2}\frac{1}{m_{2s-1}^2}\frac{1}{m_{2s-2}^2}
\frac{1}{m_{2s-3}}(1-\frac{1}{m_{2s-3}})(1-\frac{1}{m_{2s-4}}))+\cdots +\\
&& ((\prod_{j=0}^{2k}\frac{1}{m_{2s-j}})(1-\frac{1}{m_{2s}})(1-\frac{1}{m_{2s-(2k+1)}})+\\
&& \sum_{j=1}^k(\prod_{i=0}^{j-1}\frac{1}{m_{2s-i}^2})(\prod_{i=j}^{2k-j}\frac{1}{m_{2s-i}})(1-\frac{1}{m_{2s-j}})(1-\frac{1}{m_{2s-(2k+1)+j}}))+\cdots
\end{eqnarray*}

If $m_{2s}\ge 3$, then $\frac{1}{m_{2s}^2}\le \frac{1}{m_{2s}}(1-\frac{1}{m_{2s}})\le \frac{2}{9}$. Hence, the above product is at most
\begin{eqnarray*}
&& \frac{2}{9}((1-\frac{1}{m_{2s-1}})+ \frac{1}{m_{2s-1}}\frac{1}{m_{2s-2}}(1-\frac{1}{m_{2s-3}})+\frac{1}
{m_{2s-1}}(1-\frac{1}{m_{2s-1}})(1-\frac{1}{m_{2s-2}})\\
&&+\frac{1}{m_{2s-1}}\frac{1}{m_{2s-2}}\frac{1}{m_{2s-3}}\frac{1}{m_{2s-4}}(1-\frac{1}{m_{2s-5}})
+\frac{1}{m_{2s-1}}\frac{1}{m_{2s-2}}\frac{1}{m_{2s-3}}(1-\frac{1}{m_{2s-1}})(1-\frac{1}{m_{2s-4}})\\
&&+\frac{1}{m_{2s-1}^2}\frac{1}{m_{2s-2}}(1-\frac{1}{m_{2s-2}})(1-\frac{1}{m_{2s-3}})+ \frac{1}{m_{2s-1}}\frac{1}{m_{2s-2}}\frac{1}{m_{2s-3}}
\frac{1}{m_{2s-4}}\frac{1}{m_{2s-5}}\frac{1}{m_{2s-6}}(1-\frac{1}{m_{2s-7}})\\
&&+\frac{1}{m_{2s-1}}\frac{1}{m_{2s-2}}\frac{1}{m_{2s-3}}
\frac{1}{m_{2s-4}}\frac{1}{m_{2s-5}}(1-\frac{1}{m_{2s-1}})(1-\frac{1}{m_{2s-6}})\\
&&+\frac{1}{m_{2s-1}^2}\frac{1}{m_{2s-2}}\frac{1}{m_{2s-3}}
\frac{1}{m_{2s-4}}(1-\frac{1}{m_{2s-2}})(1-\frac{1}{m_{2s-5}})\\
&&+\frac{1}{m_{2s-1}^2}\frac{1}{m_{2s-2}^2}\frac{1}{m_{2s-3}}
(1-\frac{1}{m_{2s-3}})(1-\frac{1}{m_{2s-4}})+\cdots +\\
&&((\prod_{j=1}^{2k}\frac{1}{m_{2s-j}})(1-\frac{1}{m_{2s-(2k+1)}})+\sum_{j=1}^k(\prod_{i=1}^{j-1}\frac{1}{m_{2s-i}^2})(\prod_{i=j}^{2k-j}\frac{1}{m_{2s-i}})(1-\frac{1}{m_{2s-j}})(1-\frac{1}{m_{2s-(2k+1)+j}}))+\cdots \\
&&=\frac{2}{9}(1-\frac{1}{m_{2s-1}^2} +\frac{2}{m_{2s-1}^2}\frac{1}{m_{2s-2}}-\frac{1}{m_{2s-1}^2}\frac{1}{m_{2s-2}^2}-\frac{2}{m_{2s-1}^2}\frac{1}
{m_{2s-2}}\frac{1}{m_{2s-3}}+\frac{2}{m_{2s-1}^2}\frac{1}{m_{2s-2}^2}\frac{1}{m_{2s-3}}\\
&&+\frac{2}{m_{2s-1}^2}\frac{1}{m_{2s-2}}\frac{1}{m_{2s-3}}\frac{1}{m_{2s-4}}-+\cdots \\
&&-\prod_{j=1}^k\frac{1}{m_{2s-j}^2}-2\sum_{j=1}^{k-1}(\prod _{i=1}^j\frac{1}{m_{2s-i}^2})(\prod_{i=j+1}^{2k-j}\frac{1}{m_{2s-j}})+2\sum_{j=1}^k(\prod_{i=1}^j\frac{1}{m_{2s-i}^2})(\prod_{i=j+1}^{2k+1-j}\frac{1}{m_{2s-j}})-+\cdots )\le \\
&&\le \frac{2}{9}.
\end{eqnarray*}
Let us suppose that $m_{2s}=2$ if $s$ is large enough. If $m_{2t-1}\ge 3$ for infinitely many $t$, then
\begin{eqnarray*}
\liminf_{s\to \infty }&& ( \frac{1}{m_{2s}}(1-\frac{1}{m_{2s-1}}) \frac{m_{2s}}{m_{2s-1}}\frac{m_{2s-2}}{m_{2s-3}}\cdots \frac{m_2}{m_1}\\
&&+\frac{1}{m_{2s-1}^2}\frac{1}{m_{2s-2}}(1-\frac{1}{m_{2s-3}}) \frac{m_{2s-2}}{m_{2s-3}}\frac{m_{2s-4}}{m_{2s-5}}\cdots \frac{m_2}{m_1}\\
&&+\frac{1}{m_{2s-1}^2}\frac{1}{m_{2s-3}^2}\frac{1}{m_{2s-4}}(1-\frac{1}{m_{2s-5}}) \frac{m_{2s-4}}{m_{2s-5}}\frac{m_{2s-6}}{m_{2s-7}}\cdots \frac{m_2}{m_1}+\cdots )=0.
\end{eqnarray*}
\vskip2mm
Thus, we may assume that $m_{t}=2$ for $t\ge 3$. Then
\begin{eqnarray*}
\liminf_{s\to \infty }&& ( \frac{1}{m_{2s}}(1-\frac{1}{m_{2s-1}}) \frac{m_{2s}}{m_{2s-1}}\frac{m_{2s-2}}{m_{2s-3}}\cdots \frac{m_2}{m_1}\\
&&+\frac{1}{m_{2s-1}^2}\frac{1}{m_{2s-2}}(1-\frac{1}{m_{2s-3}}) \frac{m_{2s-2}}{m_{2s-3}}\frac{m_{2s-4}}{m_{2s-5}}\cdots \frac{m_2}{m_1}\\
&&+\frac{1}{m_{2s-1}^2}\frac{1}{m_{2s-3}^2}\frac{1}{m_{2s-4}}(1-\frac{1}{m_{2s-5}}) \frac{m_{2s-4}}{m_{2s-5}}\frac{m_{2s-6}}{m_{2s-7}}\cdots \frac{m_2}{m_1}+\cdots )=\frac{1}{3}\frac{m_2}{m_1}
\end{eqnarray*}
and
\begin{eqnarray*}
\liminf_{s\to \infty }&&((1-\frac{1}{m_{2s}}) \frac{m_{2s-1}}{m_{2s}}\frac{m_{2s-3}}{m_{2s-2}}\cdots \frac{m_1}{m_2}\\
&&+\frac{1}{m_{2s}^2}(1-\frac{1}{m_{2s-2}}) \frac{m_{2s-3}}{m_{2s-2}}\frac{m_{2s-5}}{m_{2s-4}}\cdots \frac{m_1}{m_2}\\
&&+\frac{1}{m_{2s}^2}\frac{1}{m_{2s-2}^2}(  1-\frac{1}{m_{2s-4}}) \frac{m_{2s-5}}{m_{2s-4}}\frac{m_{2s-7}}{m_{2s-6}}\cdots \frac{m_1}{m_2}+\cdots ) \} =\frac{2}{3}\frac{m_1}{m_2}.
\end{eqnarray*}
If $x>0$ and $y>0$, then $\min\{ x,y\}\le \sqrt{xy}$ and we are done.
\vskip2mm

This completes the proof of Theorem \ref{thm4}. \hfill$\Box$\\


\begin{thebibliography}{99}

\bibitem{Erdos} P. Erd\H os, R. Freud, On disjoint sets of differences, J. Number Theory 18 (1984), 99-109.

\bibitem{Narkiewicz} W. Narkiewicz, Remarks on a conjecture of Hanani in additive number theory, Colloq. Math. 7 (1959/60), 161-165.

\bibitem{Sandor} J.H. Fang, C. S\'andor, On sets with sum and difference structure, arXiv:2205.06553.

\end{thebibliography}
\end{document}